\documentclass{amsart}
\usepackage{amssymb,amsthm,amsmath}
\usepackage[a4paper]{geometry}
\newcommand{\M}{\mathcal{M}}
\newcommand{\N}{\mathcal{N}}
\newcommand{\RR}{\mathbb{R}}

\newcommand{\C}{\mathcal{C}}

\newcommand{\A}{\mathcal{A}}
\newcommand{\R}{\mathcal{R}}
\newcommand{\E}{\mathcal{E}}

\newcommand{\D}{\partial}

\newcommand{\<}{\langle}

\newcommand{\vol}{\mathrm{vol}}

\renewcommand{\>}{\rangle}
\let\div\relax
\DeclareMathOperator{\div}{\mathrm{div}}
\DeclareMathOperator{\grad}{\mathrm{grad}}
\DeclareMathOperator{\li}{\mathrm{li}}
\newtheorem{thm}{Theorem}
\newtheorem{lem}[thm]{Lemma}

\theoremstyle{definition}

\title{Detecting anisotropic inclusions through EIT}
\author{Jan Cristina}
\address{Section de Math\'ematiques\\
  \'Ecole Polytechnique F\'ed\'erale de Lausanne
}
\email{jan.cristina@epfl.ch}

\author{Lassi P\"aiv\"arinta}
\address{Department of Mathematics\\ Tallin University of Technology}
\email{lassi.paivarinta@ttu.ee}
\thanks{Jan Cristina was supported by the Finnish Academy Center of Excellence in Inverse Problems, and the Finnish Cultural Foundation grant ``Harmonic maps, coordinate gauges and anisotropic inverse problems''.  Lassi P\"aiv\"arinta was supported by the European Research Council-2010 Advanced Grant 267700 - InvProb (Inverse Problems), and Estonian Research Council Grant PUT1093}
\date{}
\begin{document}
\maketitle
\begin{abstract}We study the evolution equation $\D_{t}u=-\Lambda_{t}u$ where $\Lambda_{t}$ is the Dirichlet--Neumann operator of a decreasing family of Riemannian manifolds with boundary.  We derive a lower bound for the solution of such an equation, and apply it to a quantitative density estimate for the restriction of harmonic functions on $\M=\Sigma_{0}$ to the boundaries of $\D\Sigma_{t}$. Consequently we are able to derive a lower bound for the difference of the Dirichlet--Neumann maps in terms of the difference of  a background metrics $g$ and an inclusion metric  $g+\chi_{\Sigma}(h-g)$ on a manifold $\M$.
\end{abstract}
The Calder\'on problem asks whether one can determine the coefficients $\sigma$ of a boundary value problem on a Riemannian manifold $(\M,g)$ 
\begin{align}
  \div(\sigma du)=0 \qquad u|_{\D\M}=f
  \label{eqn:BVP_definition}
\end{align}
from knowledge of the map $\Lambda_{\sigma}:f\mapsto \D_{\nu}u$, where $\D_{\nu}$ is the outward unit normal.  This map is called the Dirichlet--Neumann (DN) map.  In the event that $\sigma$ is a positive definite tensor, and the dimension of the manifold is greater than 2, then this is equivalent to determining a Riemannian metric from the Dirichlet--Neumann map. 

The DN-map can be formulated weakly as
\begin{equation}
  \label{eqn:DN_weak_form}
  \int_{\D\M}\Lambda_{g}(f_{1})(f_{2})\:d\vol_{\D\M}=\int_{\M}du_{1}\wedge\star _{g}d u_{2},
\end{equation}
where $\star_{g}$ is the Hodge star associated to $g$.  It can be seen that very little regularity must be assumed, and in fact a bounded measurable metric is already enough to pose the question, as then we are guaranteed existence and uniqueness for the boundary value problem \eqref{eqn:BVP_definition}.  

We study a Riemannian manifold $(\M,g)$, and its associated Dirichlet--Neumann map $\Lambda$.  Given a smooth submanifold with boundary of the same dimension, $\Sigma$ and a Riemannian metric $h$ on $\Sigma$, we can consider $\M$ with the bounded measureable metric $g+\chi_{\Sigma}(h-g)$, where $\chi_{\Sigma}$ is the characteristic function of $\Sigma$.  The main result of the paper is the following Theorem:

\begin{thm}
  \label{thm:Riem_inc_est}
  Let $\Sigma\subset \subset \M$ be a compactly contained $C^{4}$ subdomain which is an $n$-dimensional submanifold with boundary in a compact $n$-dimensional $C^{4}$ manifold with boundary $\M$. Let  $(\M,g)$ $(\Sigma,h)$ be $C^{3}$-smooth Riemannian manifolds.  Let $g/R\leq h\leq Rg$ and let $h$ be $K$ Lipschitz (with respect to $g$) on $\Sigma$.  There are positive numbers $C_{1}$, and $C_{2}$ which depend on $\M$, $g$, $R$, $\Sigma$ and $K$ such that
  
  %Given a point $x\in \D\Sigma$, there are positive numbers $C_{1}$, $C_{2}$  and $C_{3}$ which depend on $\M,g,\Sigma_{1},R,x$ and $K$,  and functions $u,v\in H^{1/2}(\D\M)$ $\|u\|_{1/2}=\|v\|_{1/2}=1$, such that
  %\[(\Lambda_{(\M,g)}-\Lambda_{(\M,g+(h-g)\chi_{\Sigma})})(u)(v)\geq C_{1}|h-g|_{g}(x)e^{-C_{2}(|h-g|_{g}(x))^{-C^{3}}}.\]
  \[\sup_{x\in \D\Sigma}|h-g|_{g}(x)\leq C_{1}|\log\|\Lambda_{\M,g}-\Lambda_{\M,h}\||^{-1/C_{2}},\]
  where $|h-g|_{g}(x)=\sup_{X\in T_{x}\M}|h(X,X)-g(X,X)|/g(X,X).$
\end{thm}
Calder\'on first posed his problem in \cite{Calderon_Problema_classic} (republished in \cite{Calderon_Problem_new}).  There he considered the problem of determining the conductivity of an object, by measuring the resulting current arising from a voltage distribution applied to the surface.  This was modelled by a Euclidean domain $\Omega$ and a scalar conductivity, $\sigma$, so that the resulting voltages satisfy the conductivity equation $\div(\sigma\grad u)=0.$  Since then much progress has been made in the case of scalar conductivities, \cite{Sylvester_Uhlmann},\cite{Calderon_in_plane},\cite{Lipschitz_conductivities},\cite{CaroRogers_Lipschitz_conductivities}.  

Also of interest are anisotropic conductivities, where $\sigma$ is a symmetric positive definite tensor; these provide a more accurate model of certain materials e.g. muscle or nerve tissue \cite{cat_dorsal_column_anisotropy,skeletal_muscle_anisotropy}). The most complete analysis of this problem in two-dimensions is \cite{anisotropic_Calderon_plane}.  In higher dimensions, as mentioned, this is equivalent to trying to determine the Riemannian metric of a manifold with boundary.  An obvious first obstruction, of course, is the invariance of the problem under boundary fixing diffeomorphisms.  I.e. given $(\M,g)$ and a boundary fixing diffeomorphism $\Phi$, then 
\[\Lambda_{g}=\Lambda_{\Phi^{*}g}.\]
Progress on this front, has been made in the case of real-analytic manifolds and metrics \cite{Lassas_Uhlmann_Taylor}, and in distinguishing between representatives within a conformal class of certain admissible manifolds in \cite{Limiting_Carleman_weights}.  But the methods used therein require rather strong assumptions on the regularity and  geometry  of the manifold respectively. Of course the results are significantly stronger.

The proof of the main theorem follows from a quantitative density estimate for the restriction of harmonic functions to the boundary of a manifold on the interior of $\M$.  

More precisely, we construct a map $\varphi:\D\M\times [0,1]\to \M$, which is a diffeomorphism onto its image, for which $ \varphi(x,0)=x$.  Denoting by $\Sigma_{t}$ the set $\M\setminus \varphi(\D\M\times[0,t))$, and assuming that $\Sigma\subset \Sigma_{1}$ and $\D\Sigma\cap \D\Sigma_{1}$ is an $(n-1)$ manifold in a neighbourhood of a point $x$, we use the following result to prove the theorem.

\begin{thm}
  \label{thm:harmonic_density_estimate}
  Let $\varphi:[0,1]\times\D\M\to \M$ be a $C^{4}$  diffeomorphism onto its image, let $g$ be a $C^{2,1}$ metric tensor.  Let $\Sigma_{t}=\M\setminus \varphi(\D\M\times[0,t))$.  There are numbers $C_{1}$ and $C_{2}$ depending on $\varphi$, $\M$ and $g$, such that for every $f\in H^{1}(\Sigma_{1})$ $\|f\|_{1/2}=1$ and every $\varepsilon>0$ there is a $u\in H^{1/2}(\D\M)$ such that
  \[\|\E(u)|_{\Sigma_{1}}-f\|_{1/2}\leq \varepsilon\]
  and $\|u\|_{1/2}\leq e^{C_{1}(\|f\|_{1}/(\varepsilon\|f\|_{1/2}))^{C_{2}}}$.
\end{thm}

The proof of theorem \ref{thm:harmonic_density_estimate} involves several key observations.  The first is that given a harmonic function $u:\M\to \RR$ satisfying $u|_{\D\M}=f$, it satisfies the \emph{tautological evolution equation}
  \begin{equation}
    \D_{t}u_{t}=-\eta_{t}\Lambda_{t}u_{t}+X_{t}u_{t}\qquad u_{0}=f,
    \label{eqn:tautological_evln_eqn}
  \end{equation}
  where $u_{t}(x)=u(\varphi(x,t))$, $\Lambda_{t}$ is Dirichlet--Neumann operator for $(\Sigma_{t},g)$ whose unit outward normal is given by $\nu_{t}$, $\eta_{t}=g(\varphi_{*}\D_{t}, \nu_{t})$,  $X_{t}$ is a derivation given by the projection of $\varphi_{*}\D_{t}$ to $\D\Sigma_{t}$. 

  We are able to derive a lower bound for solutions of \eqref{eqn:tautological_evln_eqn}, by considering the evolution of $(\|u_{t}\|_{1}/\|u_{t}\|_{0})^{2}$ and judiciously applying Gr\"onwall's lemma.

  The density result then follows by iterating the lower bound with evolution from the exterior boundary to the interior boundary, and vice versa.

  Theorem \ref{thm:Riem_inc_est} then follows by constructing functions on $\D\Sigma_{1}$ which are supported in a neighbourhood of $x$, and see the inclusion at $\D\Sigma$.  Care must be shown when choosing the functions as the metric $g+(h-g)\chi_{\Sigma}$ is not smooth on $\Sigma_{1}$, and hence Theorem \ref{thm:harmonic_density_estimate} cannot be applied as is.  

\section{Preliminaries}
We define Sobolev spaces on our manifolds via a smooth (that is as smooth as the manifold allows) partition of unity $\psi_{i}$, supported on a set $U_{i}$ with a coordinate chart $\varphi_{i}:U_{i}\to \RR^{n}_{+}$.  A measurable function $u$ is in $H^{s}(\M)$ if $\psi_{i}\circ u\circ \varphi_{i}$ is in $H^{s}(\RR^{n}_{+})$ for every $s$.  If $\M$ is $C^{k}$ then this is valid for $s\in [0,k]$. For brevity's sake whenever $C^{k,1}$ is implied for $k=-1$ we intend it to mean $L^{\infty}$.

For $s\in [0,k-1]$ we can define the spaces $H^{s}(\M,\Lambda^{l}\M)$ to be the space of measurable $l$-forms $\alpha$ for which $\varphi_{i}^{*}\psi_{i}\alpha\in H^{s}(\RR^{n}_{+},\Lambda^{l}\RR^{n})$. 

We define the spaces $H^{s}(\D\M)$ and $H^{s}(\D\M,\Lambda^{l}\D\M)$ similarly.  We define the space 
\[H^{-s}(\D\M,\Lambda^{l}\D\M)\]
to be the dual of $H^{s}(\D\M, \Lambda^{n-1-l}\D\M)$ for $s\in[0,k-1]$.  This negates the need for a volume form which will be important in the low regularity case for $s\in[0,k-1]$.  If $l=0$ then this definition can be extended to $s\in[0,k]$.  A norm $\|\cdot\|_{s}$ is fixed, although it is not particularly important which one.  For instant $\<(I+\Delta)^{s}\cdot,\cdot\>$ where $\Delta$ is some fixed Laplace--Beltrami operator.

Let $\M$ be $C^{k+1}$, Given a metric tensor $g\in C^{k-1,1}$, let $s\in[1/2,k+3/2]$. For every $f\in H^{s}(\D\M)$ the  boundary value problem 
\begin{equation}
  \int_{\M}du\wedge\star d\varphi=0\qquad u|_{\D\M}=f
  \label{eqn:BVP_forms}
\end{equation}
has a unique solution $u$ which is in $H^{s+1/2}(\M)$.  This is a standard result in second order elliptic partial differential equations \cite{GilbargTrudinger}.  We define $\E_{g}(f):=u$ to be the solution operator. The Dirichlet--Neumann map is given by 
 \begin{equation}
   \int_{\D\M}\Lambda_{g}(f_{1})(f_{2})\:d\vol=\int_{\M}d\E_{g}(f_{1})\wedge\star_{g}d\E_{g}(f_{2}).
    \label{eqn:DN_defnition})
  \end{equation}
  If $g$ is only bounded and measurable it is only defined as a map $H^{1/2}(\D\M)\to H^{-1/2}(\D\M,\Lambda^{n-1})$ and we write $\Lambda_{g}(f_{1})(f_{2})$.

\begin{lem}
  \label{lem:DN_principle_symbol}Let $k\geq 1$.  Let $g$ be a $C^{k-1,1}$ metric tensor, and $\M$ is $C^{k+1}$. The Dirichlet--Neumann map for $(\M,g)$ is given by
  \[\Lambda=\Delta^{1/2}_{0}+R,\]
  where $R:H^{s}\to H^{s}$ for $s\in [-k-1/2,k+1/2]$,  and $\Delta_{0}$ is the Laplace--Beltrami operator on $\D\M$.
\end{lem}
\proof
We will prove that $\E(u)=\E_{0}(u)+\R(u)$ where $\E_{0}:H^{s+1/2}\to H^{s+1}$ and $ \R:H^{s+1/2}\to H^{s+2}$. First let $\varphi:\D\M\times[0,\varepsilon]$ be a $C^{k+1,1}$ diffeomorphism into $\M$, such that $\varphi_{*}\D_{t}=-\nu$ at $\D\M$.  Such a difeomorphism exists by the tubular-neighbourhood or collar neighbourhood theorems \cite{Collar_neighbourhood} \cite[Thm 6.17]{Lee_smooth_manifolds}.  Then define an extension function $\E_{0}$ by taking 
\[\E_{0}u_{\lambda}(\varphi(x,t))=\psi(t)e^{-\lambda t}u_{\lambda}(x),\]
where $\psi$ is a smooth non-negative cut-off function equal to one in a neighbourhood of $0$, and $0$ in a neighbourhood of $\varepsilon$, and $\Delta_{0} u_{\lambda}=\lambda^{2}u_{\lambda}$.

Then consider $\E(u)=\E_{0}(u)+h$ where $h\in H^{1}_{0}(\M)$.  But $\E(u)$ is the solution operator for $\Delta$, the Laplace--Beltrami operator on $\M$, so 
\[\Delta h=-\Delta \E_{0}u.\]
Because $g\in C^{k,1}$, $u_{\lambda}\in C^{k+1}(\D\M)\cap H^{k+2}(\D\M)$ \cite{GilbargTrudinger}.  Consequently, $\E_{0}(u_{\lambda})\in H^{k+2}(\M)$.  Then
\[\Delta\E_{0}(u_{\lambda})=(\Delta-\Delta_{0}+\D_{t}^{2})\E_{0}(u_{\lambda})+(-\D_{t}^{2}+\Delta_{0})(\E_{0}(u_{\lambda})),\]
Note that the leading coefficients of $(\Delta_{0}-\Delta-\D_{t}^{2})$ are bounded by a constant times $t$, as our collar neighbourhood is $C^{2}$. Note further that for every $0<s\leq k+2$, there is a $C_{s}$ such that for every $\lambda$, $\|u_{\lambda}\|_{s}\leq C_{s}\lambda^{s}$. Hence
\[\|X^{k}(\Delta_{0}-\Delta)\E_{0}(u_{\lambda})(\cdot,t)\|_{L^{2}(\D\M)}\leq C_{k}( t e^{-\lambda t}\lambda^{k+2}+e^{-\lambda t} \lambda^{k+1}).\]

We will show that $\Delta\E_{0}$ extends to a bounded map from $H^{s+1/2}(\D\M)\to H^{s}(\M)$, and consequently by inverting $\Delta$ on $H^{1}_{0}(\M)$ and applying standard elliptic theory, we get that $\Delta^{-1}(\Delta\E_{0}):H^{s+1/2}(\D\M)\to H^{s+2}(\M)$.  Now we can take the normal derivative to get 
\[\D_{\nu}\Delta^{-1}(\Delta\E_{0}):H^{s+1/2}(\D\M)\to H^{s+1}(\M),\]
finally taking the trace yields a map $H^{s+1/2}(\D\M)\to  H^{s+1/2}(\D\M)$. Lastly because $\Delta^{1/2}$ and $\Lambda$ are selft adjoint it follows that $R$ is self adjoint and hence bounded for $s\in[-k-1/2,k+1/2]$.

\endproof
\begin{lem}
  \label{lem:DN_definition}  Let $k\geq 0$. Let $\M$ be a $C^{k+2}$ manifold.  Let $g\in C^{k-1,1}$ be a metric tensor on $\M$, then the Dirichlet--Neumann map $\Lambda_{g}:H^{s+1}(\D\M)\to H^{s}(\D\M,\Lambda^{n-1}\D\M)$ is  a bounded invertible map $H^{s}(\D\M)/\C\to H^{s-1}/\C$ for $s\in[1/2,k+1/2]$, where $\C$ is the space of constant functions.
\end{lem}
\proof 
We note that 
\[\Lambda_{g}(f)(f)=\int_{\M}d\E(f)\wedge\star_{g}d\E(f)=\|d\E(f)\|_{2}^{2}\]
Let $u=\E(f)$, and let $u_{\M}$ denote its average over $\M$.  By the Poincar\'e inequality for $\M$ we have
\[\|u-u_{\M}\|^{2}\leq C\int_{\M}du\wedge\star_{g}du=C\Lambda_{g}(f)(f).\]
Furthermore letting $\eta$ denote the solution to the Dirichlet problem
\[\Delta \eta= (\int_{\M}d\vol_{\M})^{-1}\qquad \eta|_{\D|M}=0,\]
we see that $u_{\M}=\int_{\D\M}f i^{*}\star_{g} d\eta$, and hence
\[\|u\|_{1}^{2}\leq C[\|f\|_{0}^{2}+\Lambda(f)(f)].\]
Then noting that the trace operator is bounded from $H^{1}(\M)\to H^{1/2}(\D\M)$ yields 
\[\|f\|^{2}_{1/2}\leq C[\|f\|_{-1/2}^{2}+\Lambda(f)(f)]\]
\[\|f\|_{1/2}\leq C[\|f\|_{-1/2}+\|\Lambda(f)\|_{-1/2}]\]
Consequently $\Lambda_{g}$ has a discrete spectrum.  The kernel of $\Lambda_{g}$ is given by the constant functions, otherwise the double extension on the doubled manifold is a non-constant harmonic function.  Hence there is a spectral gap and $\Lambda$ is invertible on the orthogonal complement of the constant functions. If $k=0$ then we are done.

For $k\geq 1$ consider $[\Lambda,\Delta_{0}^{\sigma/2}]$.  This is equal to $[R,\Delta_{0}^{\sigma/2}]$, which is a bounded map $H^{s}\to H^{s-\sigma}$ for $s\in[0,k+1/2]$, and $\sigma/2\in [0,k+1/2]$.  So we can apply the preceding to yield
\begin{align*}
  \|\Delta_{0}^{\sigma/2}f\|_{1/2}&\leq C\|\Delta_{0}^{\sigma/2} f\|_{-1/2}+\|\Delta_{0}^{\sigma/2}\Lambda(f)\|_{1/2}+\|[\Lambda,\Delta_{0}^{\sigma/2}] f\|_{-1/2}\\
  &\leq C[\|f\|_{\sigma-1/2}+\|\Lambda(f)\|_{\sigma-1/2}]
\end{align*}
As a result
\[\|f\|_{s+1/2}\leq C[\|f\|_{s-1/2}+\|\Lambda(f)\|_{s-1/2}]\]
for $s\in [0,k-1/2]$.
Consequently we can deduce that the eigenfunctions of $\Lambda\in H^{s}(\D\M)$ for $s\in [0,k+1/2]$, and $\Lambda$ has a bounded inverse $H^{s-1}\to H^{s}$.
\endproof

\begin{lem}
   Let $k\geq 0$ and $g\in C^{k,1}$, $\eta\in C^{k+1,1}(\D\M)$, then $[\Lambda,\eta]:H^{s}\to H^{s}$ for  $s\in[0,k+1]$
  \label{lem:DN_scalar_commutator}
\end{lem}
This is a classical result in pseudodifferential operators \cite{psidos_NLPDE} but results with less than $C^{\infty}$ regularity are hard to find, hence we prove it here.
\proof
The bulk of the proof is to show that $\Delta_{0}^{1/2}\eta-\eta\Delta_{0}^{1/2}$ is bounded from $H^{s}\to H^{s}$.  This relies on the fact that Dirichlet Neumann map can be expressed as
\begin{multline}
  \int_{\D\M}\Delta_{0}^{1/2}(u)v\:d\vol=\int_{\M}d\E(u)\wedge d\E(v)\\ =\int_{\M}d(\E(u)+h)\wedge \star d\E(v)=\int_{\M}d\E(u)\wedge\star d(\E(v)+h')
\end{multline}
where $\E$ is the harmonic extension operator, and $h$ and $h'$ are arbitrary $W^{1,2}_{0}$ functions.  Thus
\begin{align*}
  \int_{\N}\Delta_{0}^{1/2}(\eta u_{\lambda})u_{\mu}&-\eta\Delta_{0}^{1/2}(u_{\lambda})u_{\mu}\:d\vol\\
  &=\int_{\N}\eta u_{\lambda}\Delta_{0}^{1/2}(u_{\mu})-\Delta_{0}^{1/2}(u_{\lambda})\eta u_{\mu}\\
  &=\int_{\N}\int_{0}^{\infty}d(\eta e^{-t\lambda}u_{\lambda})\wedge\star d(e^{-t\mu u_{\mu}})-d(e^{-t\lambda}u_{\lambda})\wedge\star d(\eta e^{-t\mu}u_{\mu})\\
  &=\int_{N}\int_{0}^{\infty}d\eta\wedge\star(u_{\lambda}du_{\mu}-u_{\mu}du_{\lambda})e^{-t(\mu+\lambda)}\:dt\\
  &\leq C \|d\eta\|_{\infty}(\|u_{\lambda}\|_{1}+\|u_{\mu}\|_{1})(\lambda+\mu)^{-1}.
\end{align*}
 To show that it is bounded from $H^{s}\to H^{s}$, we consider only integer $s$, as fractional $s$ follows from interpolation.  We can commute with  $\Delta^{s/2}$.  If $s$ is odd this is
\begin{align*}
  \Delta_{0}^{s/2}[\Delta_{0}^{1/2},\eta]&= \eta\Delta_{0}^{(s+1)/2}+[\Delta_{0}^{(s+1)/2},\eta]-\Delta_{0}^{s/2}\eta \Delta_{0}^{1/2}\\
  &=[\eta,\Delta_{0}^{1/2}]\Delta_{0}^{s/2}+[\Delta_{0}^{(s+1)}/2,\eta]+\Delta_{0}^{1/2}[\Delta_{0}^{(s-1)/2},\eta]\Delta_{0}^{1/2}\\
  &=[\eta,\Delta_{0}^{1/2}]\Delta_{0}^{s/2}+B
\end{align*}
where $B$ is bounded $H^{s}\to H^{0}$.
If $s$ is even, then
\begin{align*}
  \Delta_{0}^{s/2}[\Delta_{0}^{1/2},\eta]&=\Delta_{0}^{1/2}\eta\Delta_{0}^{s/2}+\Delta_{0}^{1/2}[\Delta_{0}^{s/2},\eta]-\eta\Delta_{0}^{1/2}\Delta_{0}^{s/2}+[\Delta_{0}^{s/2},\eta]\Delta_{0}^{1/2}\\
  &=[\Delta_{0}^{1/2},\eta] \Delta_{0}^{s/2}+B
\end{align*}
where $B:H^{s}\to H^{0}$. 
\endproof
\begin{lem}
  For $k\geq 0$ let $g$ be $C^{k,1}$. Let $X$ be a derivation on $\D\M$, with $C^{k,1}$ coefficients.  Then 
  \[ [X,\Lambda]:H^{s}\to H^{s-1},\]
 For $s\in[1/2,k]$.
  \label{lem:DN_vector_commutator}
\end{lem}
\proof
Once again
\begin{align*}
  \int_{\D\M}[X,\Delta_{0}^{1/2}](u)(v)\:d\vol&=-\int_{\N}\Delta_{0}^{1/2}(u)X(v)\:d\vol+\int_{\N}\Delta_{0}^{1/2}(u)\div X v\:d\vol \\
  &\quad-\int_{\N} X(u)\Delta_{0}^{1/2}(v)\:d\vol
\end{align*}
Then we examine the action on basis elements, extending $X(u_{\lambda})$ by $e^{-\lambda t}X(u_{\lambda})$.  Consequently, using an orthonormal frame $e_{i}$,
\begin{align*}
  \int_{\D\M}[X,&\Delta_{0}^{1/2}](u_{\lambda})(u_{\mu})\:d\vol\\
  &=-\int_{0}^{\infty}\int_{\D\M} \lambda \mu e^{-t(\lambda+\mu)}[X(u_{\lambda})u_{\mu}+u_{\lambda}(X-\div X)(u_{\mu})]\:d\vol\:dt\\
  &\quad-\int_{0}^{\infty}e^{-t(\lambda+\mu)}\int_{\D\M}[d(X(u_{\lambda}))\wedge\star d(u_{\mu})+d(u_{\lambda})\wedge\star d((X-\div(X))u_{\mu})]\:dt\\
  &=\frac{1}{\lambda+\mu}\int_{\D\M}[e_{i}(X(u_{\lambda}))e_{i}(u_{\mu})+e_{i}(u_{\lambda})e_{i}( (X-\div X)(u_{\mu}))]\:d\vol\\
  &=\frac{1}{\mu+\lambda}\int_{\D\M}e_{i}(X(u_{\lambda}))e_{i}(u_{\mu})+e_{i}(u_{\lambda})[X(e_{i})(u_{\mu})+[e_{i},X](u_{\mu})\\
&\quad-e_{i}(\div X) u_{\mu}-\div X e_{i}(u_{\mu})]\:d\vol\\
&=\frac{1}{\mu+\lambda}\int_{\D\M}[e_{i},X](u_{\lambda})e_{i}(u_{\mu})+e_{i}(u_{\lambda})[e_{i},X](u_{\mu})-e_{i}(u_{\lambda})u_{\mu}e_{i}(\div X)\:d\vol.
\end{align*}

So $[X,\Delta_{0}^{1/2}]$ is bounded from $H^{1/2}\to H^{-1/2}$.  Suppose $k$ is even and consider that
\[[\Delta_{0}^{k/2},[X,\Delta_{0}^{1/2}]]=-[\Delta_{0}^{1/2},[\Delta_{0}^{k/2},X]],\]
but because $X\in C^{k,1}$ $[\Delta_{0}^{k},X]$ is a differential operator of order $k$, and so is bounded from $H^{k-1/2}\to H^{-1/2}$.  Odd $k$ follows as in the case of the proof of Lemma \ref{lem:DN_scalar_commutator}.
\endproof

\begin{lem}
  \label{lem:DN-derivative}  Assume $g\in C^{k,1}$  and let $\Phi:[0,t_{0}]\times\D\M\to \M$ be a $C^{k+2}$ map which is a diffeomorphism onto its image. Let $\Lambda_{t}$ denote the Dirichlet--Neumann map for $\Sigma_{t}=\M\setminus \Phi([0,t]\times\D\M)$, and let $A_{t}=\varphi^{*}_{t}\Lambda_{t}$.  Then $\D_{t}A_{t}$ exists in the weak sense
  \[\lim_{h\to 0}\frac{1}{h}((A_{t+h}-A_{t})u,v)_{t}=(\D_{t}A_{t}u,v)\]
  for every $u,v\in H^{1/2}$.  Furthermore  $\D_{t}A_{t}$ is a bounded operator $H^{s+1}\to H^{s}$ for $s\in [1/2,k+1/2]$, uniformly in $t$.
\end{lem}
\proof
The proof of this is adapted from a similar proof in \cite{Sylvester_layer_stripping}. Let $\E_{t}(f)$ denote the extension operator at $\Sigma_{t}$.  Then consider
\begin{align*}
  \int_{\D\M}\Lambda_{t+h}(u)(v)\:d\vol_{t+h}&-\Lambda_{t}(u)(v)\:d\vol_{t}\\
  &=\int_{\Sigma_{t+h}}d\E_{t+h}(u)\wedge\star_{g}d\E_{t+h}(v)-\int_{\Sigma_{t}}d\E_{t}(u)\wedge\star_{g}d\E_{t}(v)\\
  &=\int_{\Sigma_{t+h}}\bigg[d\E_{t+h}(u-\E_{t}(u)|_{\Sigma_{t+h}})\wedge\star_{g}d\E_{t+h}(v)\\
  &\qquad\qquad\qquad+d\E_{t}(u)\wedge\star_{g}d\E_{t+h}(v-\E_{t}(v)|_{\Sigma_{t+h}})\bigg]\\
  &\quad-\int_{\Sigma_{t}\setminus \Sigma_{t+h}}d\E_{t}(u)\wedge\star_{g}d\E_{t}(v).
\end{align*}
If we divide by $h$ and let $h$ tend to $0$. we arrive at 
\begin{align*}
  \int_{\D\M}(\D_{t}\Lambda_{t}(u)(v))\:d\vol_{t}&-\int_{\D\M}\Lambda_{t}(u)(v)\gamma_{t}\:d\vol_{t}\\
  &=\int_{\D\M}\Lambda_{t}(-\D_{\tau}(\E_{t}(u)(\varphi_{\tau}(x))(v)+\Lambda_{t}(u)(-\D_{\tau}\E_{t}(v)(\varphi_{\tau}(x))\:d\vol_{t}\\
  &\quad-\int_{\D\M}\eta_{t}\Delta_{t}(u)(v)+g(d\eta_{t},du)v+\eta_{t}\Lambda_{t}(u)\Lambda_{t}(v)\:d\vol_{t},
\end{align*}
  where $\gamma_{t}d\vol_{t}=\D_{t}d\vol_{t}$.
We note that $\D_{\tau}\E_{t}(u)(\varphi_{\tau}(x))=-\eta_{t}\Lambda_{t}(u)-X_{t}(u)$, and $\Lambda_{t}$ is a self adjoint operator for $d\vol_{t}$:
  \begin{align*}
    \D_{t}\Lambda_{t}=\eta_{t}(\Lambda_{t}^{2}-\Delta_{t})+[\Lambda_{t},\eta_{t}]\Lambda_{t}+[\Lambda_{t},X_{t}]+\div X_{t}\Lambda_{t}+\gamma_{t}\Lambda_{t}.
  \end{align*}
  Thus we need only apply Lemmas \ref{lem:DN_principle_symbol}, \ref{lem:DN_scalar_commutator} and \ref{lem:DN_vector_commutator}  to yield the result.
\endproof

\begin{lem}
Let $g_{1}$ and $g_{2}$ be two bounded measurable Riemannian metrics, and let $\Lambda^{1}$ and $\Lambda^{2}$ denote their respective Dirichlet--Neumann maps.  Then
\[(\Lambda^{1}-\Lambda^{2})(u)(v)=\int_{\M\setminus\Sigma_{1}}d\E^{1}(u)\wedge(\star^{1}-\star^{2})d\E^{2}(v)+(\Lambda_{\Sigma_{1}}^{1}-\Lambda_{\Sigma_{1}}^{2})(\E^{1}(u)|_{\Sigma_{1}})(\E^{2}(v)|_{\Sigma_{1}}),\]
  where $\Lambda^{i}_{\Sigma_{1}}$ is the Dirichlet--Neumann map for $g_{i}$ restricted to the submanfold $\Sigma_{1}$.
\end{lem}
\proof
First 
\[\Lambda^{i}(u)(v)=\int_{\M}d\E^{i}(u)\wedge\star^{1}d\E^{i}(v),\]
but because $\E^{1}(u)$ is $\star^{1}$ harmonic for any $h\in H^{1}_{0}(\M)$, 
\[\Lambda^{1}(u)(v)=\int_{\M}d\E^{1}(u)\wedge\star^{1}d(\E^{1}(v)+h).\]
As such because $\E^{2}(v)-\E^{1}(v)\in H^{1}_{0}$, we arrive at
\[\Lambda^{1}(u)(v)=\int_{\M}d\E^{1}(u)\wedge\star^{1}d\E^{2}(v).\]
Similarly
\[\Lambda^{2}(u)(v)=\int_{\M}d\E^{1}(u)\wedge\star^{2}d\E^{2}(v).\]
When we take the difference we arrive at
\begin{align*}
  (\Lambda^{1}-\Lambda^{2})(u)(v)&=\int_{\M}d\E^{1}(u)\wedge(\star^{1}-\star^{2})d\E^{2}(v)\\
  &=\int_{\M\setminus \Sigma_{1}}d\E^{1}(u)\wedge(\star^{1}-\star^{2})d\E^{2}(v)+\int_{\Sigma_{1}}d\E^{1}(u)\wedge(\star^{1}-\star^{2})d\E^{2}(v).
\end{align*}
The result follows from realising that the harmonic extension with respect to $g$ of $\E^{1}(u)|_{\Sigma_{1}}$ is equal to $\E^{1}(u)$ on $\Sigma_{1}$, and likewise for $\E^{2}$.
\endproof
\begin{lem}
\label{lem:extension_comparison}
Let $g$ be a $C^{0,1}$ metric on $\M$ and let $h$ be a bounded measurable metric on $\Sigma_{1}$. Let $\E^{1}:H^{1/2}(\D\M)\to H^{1}(\M)$ denote the harmonic extension operator for $g$ and let $\E^{2}$ denote that for $g+\chi_{\Sigma_{1}}(h-g)$.  Let $\E^{0,1}_{\M\setminus \Sigma_{1}}$ denote the harmonic extension operator on $\M\setminus \Sigma_{1}$ with prescribed values on $\D\Sigma_{1}$ and zero on $\D\M$, and let $\Lambda_{\M\setminus \Sigma_{1}}^{0,1}$ denote the corresponding Dirichlet--Neumann map $H^{1/2}(\D\Sigma)\to H^{-1/2}(\D\Sigma)$, then
  \[\E^{2}(v)|_{\Sigma_{1}}=(\Lambda^{0,1}_{\M\setminus\Sigma_{1}}+\Lambda_{\Sigma_{1,h}})^{-1}(\Lambda_{\Sigma_{1},g}+\Lambda^{0,1}_{g})(\E^{1}(u)|_{\Sigma_{1}})\]
  where $\Lambda_{\Sigma_{1},g}$ is the Dirichlet--Neumann map for $g$ restricted to the submanfold $\Sigma_{1}$ (respectively for $h$).
\end{lem}
\proof Let $\E^{0,1}_{g}: H^{1/2}(\D\Sigma_{1})\to H^{1}(\M\setminus \Sigma_{1})$ denote the map which takes $u$ to the harmonic function, with respect to the metric tensor $g$, which is equal to $u$ on $\D\Sigma_{1}$ and is equal to $0$ on $\D\M$.  Similarly let $\E^{1,0}_{g}:H^{1.2}(\D\M)\to H^{1}(\M\setminus \Sigma_{1})$ denote the harmonic extension operator with respect to the metric tensor $g$, which takes $u$ to the harmonic function on $\M\setminus \Sigma_{1}$ equal to $u$ on $\D\M$ and equal to $0$ on $d\Sigma_{1}$.  Let $\E_{\Sigma_{1},g},\E_{\Sigma_{1},h}:H^{1/2}(\D\Sigma_{1})\to H^{1}(\Sigma_{1}$, denote the harmonic extension operators  on $\Sigma_{1}$ with respect to the metric tensors $g$ and $h$ respectively.  We identify $\E^{1,0}_{g}(u)$ with its extension by $0$ to $\Sigma_{1}$.  And for any function $v\in H^{1/2}(\Sigma_{1})$ we denote $\E^{0,1}_{g}(v)+\E_{\Sigma_{1},g}(v)$ to be the function given by $\E^{0,1}_{g}(v)$ in $\M\setminus \Sigma_{1}$ and $\E_{\Sigma_{1},g}(v)$ in $\Sigma_{1}$.  similarly for $\E^{0,1}_{g}(v)+\E_{\Sigma_{1},h}$.

It follows that 
\[\E^{1}(u)=\E^{1,0}_{g}(u)+\E^{0,1}_{g}(u')+\E_{\Sigma_1,g}(u')\]
where $u'=\E^{1}(u)|_{\Sigma_{1}}$,
and similarly
\[\E^{2}(v)=\E^{1,0}_{g}(v)+\E^{0,1}_{g}(v')+\E_{\Sigma_{1},h}(v').\]
For $\E^{1}$ the first derivative of $\E^{1}(u)$ are continuous, so denoting $\nu_{\Sigma_{1},g}$ the outward normal of $\Sigma_{1}$ with respect to $g$, we have
\[\D_{\nu_{\Sigma_{1},g}}[\E^{1,0}_{g}(u)+\E^{0,1}_{g}(u')]=\D_{\nu_{\Sigma_{1},g}}\E_{\Sigma_{1},g}(u').\]
While $\nu_{\Sigma,g}$ is outward normal of $\Sigma_{1}$ it is the inward normal of $\M\setminus \Sigma_{1}$ at $\D\Sigma_{1}$, so
\[\D_{\nu_{\Sigma_{1},g}}(\E^{0,1}_{g}(u'))=-\Lambda^{0,1}_{g}(u'),\]
while
\[\D_{\nu_{\Sigma_{1},g}}(\E_{\Sigma_{1},g}(u'))=\Lambda_{\Sigma_{1},g}(u').\]
Putting these together we arrive at 
\begin{equation}
  \label{eqn:extension_comparison}
(\Lambda^{0,1}_{g}+\Lambda_{\Sigma_{1},g})(u')=-\D_{\nu_{\Sigma_{1},g}}\E^{1,0}_{g}(u)
.\end{equation}
Because $\Lambda^{0,1}_{g}$ and $\Lambda_{\Sigma_{1},g}$ are both Dirichlet--Neumann operators it follows that for any equivalent $H^{1/2}$ norm, there are constants $C_{1}$ and $C_{2}$ such that
\[\Lambda^{0,1}_{g}(u')(u')\geq C_{1}\|u\|_{1/2}^{2}\text{ and }\Lambda_{\Sigma_{1},g}(u')(u')\geq 0,\]
as such 
\[(\Lambda^{0,1}_{g}+\Lambda_{\Sigma_{1},g})(u')(u')\geq (C_{1}+C_{2})\|u'\|_{1/2},\]
And so $\Lambda^{0,1}_{g}+\Lambda_{\Sigma_{1},g}:H^{1/2}(\D\Sigma_{1})\to H^{-1/2}(\D\Sigma_{1})$ is invertible. 

For $\E^{2}$ we note the jump discontinuity at $\D\Sigma_{1}$ in the coefficients slightly complicates things.  However a simple inspection of the weak form of the equation yields the continuity condition
\[i^{*}\star^{1}d(\E^{1,0}_{g}(v)+\E^{0,1}_{g}(v'))=i^{*}\star^{2}d\E_{\Sigma_{1},h}(v')\]
where $i^{*}$ is the restriction of forms to $\D\Sigma_{1}$ (i.e. the pullback of the inclusion map $i:\D\Sigma_{1}\to \M$).
Consequently 
\[i^{*}\star^{1}d\E^{1,0}_{g}(v)-\Lambda^{0,1}_{g}(v')=\Lambda_{\Sigma_{1},h}(v'),\]
and so
\[i^{*}\star^{1}d\E^{1,0}_{g}(v)=(\Lambda^{0,1}_{g}+\Lambda_{\Sigma_{1},h})(v').\]
If we set $u=v$, then it follows that
\[\E^{2}(u)|_{\D\Sigma_{1}}=(\Lambda^{0,1}_{g}+\Lambda_{\Sigma_{1},h})^{-1}(\Lambda^{0,1}_{g}+\Lambda_{\Sigma_{1},g})(\E^{1}(u)|_{\D\Sigma_{1}})\]

\endproof

\begin{lem}\label{lem:extension_restriction_adjoint}Let $g$ be $C^{0,1}$, and let $\E^{1}|_{\D\Sigma_{1}}:L^{2}(\D\M,d\vol_{\D\M})\to L^{2}(\D\Sigma_{1},d\vol_{d\Sigma_{1}})$.  It is densely defined and 
\[\E^{1}|_{\D\Sigma_{1}}^{*}=\D_{\nu_{\D\M},g}\E^{0,1}_{g}(\Lambda^{0,1}_{g}+\Lambda_{\Sigma_{1},g})^{-1}\]
where $\D_{\nu\M,g}\E^{0,1}_{\M\setminus \Sigma_{1}}$ is the outward normal derivative at $\D\M$ of $\E^{0,1}_{\M\setminus \Sigma_{1}}$. 
\end{lem}\proof  Recall \eqref{eqn:extension_comparison} from the proof of Lemma \ref{lem:extension_comparison}.  This implies
\[\E^{1}(u)|_{\D\Sigma_{1}}=(\Lambda^{0,1}_{g}+\Lambda_{\Sigma_{1},g})^{-1}(\D_{\nu_{\Sigma_{1}}}\E^{1,0}_{g}(u).\]
Consequently
\[\E^{1}|_{\D\Sigma_{1}}^{*}=\D_{\nu_{\D\M},g}\E^{0,1}_{g}(\Lambda^{0,1}_{g}+\Lambda_{\Sigma_{1},g})^{-1}.\]\endproof

%We extend $\M$ by a cylinder $\D\M\times(-\infty,0]$, with the metric from extended fermi coordinates.  Again applying the reasoning for \eqref{eqn:extension_comparison}, we get
%\[\D_{\nu_{\D\M,g}}\E^{0,1}_{g}=(\Delta^{1/2}_{g_{0}}+\Lambda^{1,0}_{g}) \E^{1,0}_{hat{\M},g}.\]

%\begin{thm}
%  Let $\lambda(u)$  Then $(\E^{1}|_{\D\Sigma_{1}})^{*}(u)\geq C e^{-tC\lambda(u)}\|u\|_{1/2}$.
%\end{thm}

\section{The tautological evolution equation}
Given $\Phi:\D\M\times[0,1]\to \M$ a $C^{2}$ diffeomorphism onto its image, $\Phi(x,t)=\varphi_{t}(x)$, let $u_{t}=\E^{1}(u)\circ\varphi_{t}$, then 
\[\D_{t}u_{t}=(\varphi_{t})_{*}\D_{t}(\E^{1}(u))= (\Pi_{\nu}((\varphi_{t})_{*}\D_{t})+\Pi_{\D\Sigma_{t}}(\varphi_{t})_{*}\D_{t}) \E^{1}(u))\]
\[\D_{t}u_{t}=-\eta_{t}A_{t}u_{t}+X_{t}u_{t},\]
where $\eta_{t}=-\varphi_{t}^{*}g((\varphi_{t})_{*}\D_{t},\nu)$ and $X_{t}=(\varphi^{-1}_{t})_{*}\Pi_{\Sigma_{t}}(\phi_{t})_{*}\D_{t}$.  By the transversality of $\D_{t}$ to $\D\Sigma_{t}$, $\eta_{t}\geq \varepsilon>0$ for every $t$ and $x$.  Furthermore $\eta$ is $C^{0,1}$, and bounded away from $0$, so $\eta^{1/2}$ is also $C^{1,1}$.  $A_{t}$ is the Dirichlet--Neumann operator for $\Sigma_{t}$ conjugated with $\varphi_{t}$.  It is self adjoint with respect to $\varphi_{t}^{*}d\vol_{g,\D\Sigma_{t}}$, and bounded and invertible from $H^{s+1}\to H^{s}$ for $s\in [1/2,5/2]$.  This observation was noted in \cite{Sylvester_layer_stripping}.
\begin{align*}
  \D_{t}u_{t}=-\eta_{t}^{1/2}A_{t}\eta_{t}^{1/2}u_{t}+\eta_{t}^{1/2}[\eta^{1/2}_{t},A_{t}]u_{t}+X_{t}u_{t}.
\end{align*}
Let $B_{t}=\eta_{t}^{1/2}A_{t}\eta_{t}^{1/2}$, $S_{t}=\eta_{t}^{1/2}[\eta_{t}^{1/2},A_{t}]$
By Lemma \ref{lem:DN-derivative} $\D_{t}B_{t}=\dot{\eta}_{t} \eta_{t}^{-1/2}A_{t}\eta_{t}^{1/2}+\eta^{1/2}_{t}\dot{A}_{t}\eta^{1/2}_{t}+\eta^{1/2}_{t}A_{t}\eta^{-1/2}\dot{\eta}_{t}$.  Consequently it is bounded $H^{s+1}\to H^{s}$ for $s\in[1/2,3/2]$.

Because $\eta_{t}\in C^{1,1}$ multiplication by $\eta$ is a bounded map $H^{s}\to H^{s}$ for $s\in[0,2]$.  In particular by Lemma \ref{lem:DN_definition}  
\[\<B_{t}u_{t},B_{t}u_{t}\>=\|B_{t}u_{t}\|_{0}^{2}\sim\|u_{t}\|_{1}^{2}-\|u_{t}\|_{0}^{2}.\]
and
\[\<B_{t}u_{t},u_{t}\>=\<A_{t}\eta^{1/2}_{t}u_{t},\eta^{1/2}_{t}u_{t}\>\sim\|\eta^{1/2}_{t}u_{t}\|_{1/2}^{2}-\|\eta_{t} u_{t}\|_{0}^{2}\sim\|u_{t}\|_{1/2}^{2}-\|u_{t}\|_{0}^{2}.\]
We note that the inner product $\<\cdot,\cdot\>$ for functions on $\D\M$ obviously depends on the choice of measure, and that 
\[\D_{t}\<\cdot,\cdot\>_{t}=\<\gamma_{t}\cdot,\cdot\>.\]
This can all be formally codified into the following lemma
\begin{lem}
  Let $\Phi:\D\M\times[1,1]\to \M$ be $C^{3}$  Let $g$ be $C^{1,1}$.  For any $f\in H^{1/2}(\D\M)$, let $f_{t}=\phi_{t}\E_{1}(f)$.  Then $f\in H^{1}([0,1],H^{1/2}(\D\M),H^{1}(\D\M),L^{2}(\D\M))$ and  satisfies
  \[\D_{t}f(t)=-B_{t}f_{t}+X_{t}f_{t}+S_{t},\]
  where
  $B$ satisfies
  \[\<B_{t}u,B_{t}u\>_{t}\sim \|u\|_{1}^{2}-C\|u\|_{0}^{2}\]
  \[\<B_{t}^{*_{t}}u,v\>_{t}=\<u,B_{t}v\>_{t}\]
  \[\<B_{t}u,u\>\sim \|u\|_{1/2}^{2}-C\|u\|_{0}^{2}\]
  \[\dot{B}_{t}:H^{1}\to H^{0}\]
  \[X_{t}\text{ is a derivation with uniformly Lipschitz coefficients.}\]
  \[S_{t}:H^{1}\to H^{1},\]
all uniformly in $t$.
\end{lem}
\proof The only point not mentioned is that $f_{t}\in H^{1}([0,1],H^{1/2}(\D\M),H^{1}(\D\M),L^{2}(\D\M))$, \emph{cf.} \cite[\S 23.6]{nonlinear_fnl_analysis_and_apps_IIa} for the formal definition.  However, by the Lipschitz continuity of $g$, $f\in H^{2}(\M)$ and so $\D_{t}f\in H^{1}(\M)$ and $\D_{t}f_{t}\in L^{\infty}([0,1],H^{1/2}(\D\M))$.\endproof

\begin{lem}Let $f\in H^{1}([0,1],H^{1/2}(\D\M),H^{1}(\D\M),L^{2}(\D\M))$.  Suppose 
  \label{lem:osc_evolution}
  \begin{equation}
    \label{eqn:generic_evolution_equation}
    \D_{t}f_{t}=-B_{t}f_{t}+X_{t}f_{t}+S_{t}f_{t}
  \end{equation}
where  $B$ satisfies
  \[\<B_{t}u,B_{t}u\>_{t}\sim \|u\|_{1}^{2}-C\|u\|_{0}^{2}\]
  \[\<B_{t}^{*_{t}}u,v\>_{t}=\<u,B_{t}v\>_{t}\]
  \[\<B_{t}u,u\>\sim \|u\|_{1/2}^{2}-C\|u\|_{0}^{2}\]
  \[\dot{B}_{t}:H^{1}\to H^{0}\]
  \[X_{t}\text{ is a derivation with uniformly Lipschitz coefficients.}\]
  \[S_{t}:H^{1}\to H^{1},\]
all uniformly in $t$.  Then
\[\frac{\<B_{t}f_{t},B_{t}f_{t}\>_{t}}{\<f_{t},f_{t}\>_{t}}\leq e^{C_{1}t}\left[\frac{\<B_{0}f_{0},B_{0}f_{0}\>_{0}}{\<f_{0},f_{0}\>_{0}}+\frac{C_{2}}{C_{1}}\right]-\frac{C_{2}}{C_{1}}.\]

\end{lem}
\proof The proof follows standard methods in monotone evolution equations \cite[Theorem 23.A]{nonlinear_fnl_analysis_and_apps_IIa}. Consider
\[\lambda(t,f_{0})=\frac{\<B_{t}f_{t},B_{t}f_{t}\>_{t}}{\<f_{t},f_{t}\>_{t}}.\]
Then
\begin{align*}
  \D_{t}&\lambda(t,f_{0})\\
  &=\frac{1}{\<f_{t},f_{t}\>_{t}}\big[2\<\dot{B}_{t}f_{t},B_{t}f_{t}\>_{t}-2\<B_{t}^{2}f_{t},B_{t}f_{t}\>_{t}+2\<B_{t}X_{t}f_{t},B_{t}f_{t}\>_{t}\\
  &\quad+2\<B_{t}S_{t}f_{t},B_{t}f_{t}\>_{t}+\<B_{t}f_{t},\gamma_{t}B_{t}f_{t}\>_{t}\big]\\
  &\quad-\frac{\<B_{t}f_{t},B_{t}f_{t}\>_{t}}{\<f_{t},f_{t}\>_{t}^{2}}\big[\<B_{t}f_{t},\gamma_{t}f_{t}\>_{t}-2\<B_{t}f_{t},f_{t}\>_{t}+2\<X_{t}f_{t},f_{t}\>_{t}+2\<S_{t}f_{t},f_{t}\>_{t}\big]\\
  &=\frac{2}{\<f_{t},f_{t}\>^{2}_{t}}\overbrace{\big[-\<B^{2}_{t}f_{t}B_{t}f_{t}\>_{t}\<f_{t},f_{t}\>_{t}+\<B_{t}f_{t},B_{t}f_{t}\>_{t}\<B_{t}f_{t},f_{t}\>\big]}^{I_{1}}\\
  &\quad+\frac{1}{\<f_{t},f_{t}\>_{t}}\overbrace{\big[2\<\dot{B}_{t}f_{t},B_{t}f_{t}\>_{t}+\<B_{t}f_{t},\gamma_{t}B_{t}f_{t}\>_{t}+2\<B_{t}X_{t}f_{t},B_{t}f_{t}\>_{t}+2\<B_{t}S_{t}f_{t},B_{t}f_{t}\>_{t}\big]}^{I_{2}}\\
  &\quad+\frac{\lambda(t,f_{0})}{\<f_{t},f_{t}\>}\underbrace{\big[\<f_{t},f_{t}\>_{t}+\<f_{t},\gamma_{t}f_{t}\>_{t}+2\<\div_{t} X_{t}f_{t},f_{t}\>_{t}+2\<S_{t}f_{t},f_{t}\>_{t}\big]}_{I_{3}}.
  \end{align*}
  We claim $I_{1}\leq 0$, $|I_{2}|\leq C_{1}(\<B_{t}f_{t},B_{t}f_{t}\>_{t}+\<f_{t},f_{t}\>)$ and $|I_{3}|\leq C_{2} \<f_{t},f_{t}\>_{t}$.  For notational efficiency we temporarily drop any mention of the parameter $t$.  
  
 First we show that $I_1\leq 0$.   The operator $B$ is a self-adjoint operator for $\<\cdot,\cdot\>$ with compact resolvent, consequently it induces an eigenspace decomposition of $L^{2}$ into eigenfunctions $u_{i}$ with eigenvalues $\lambda_{i}$  Let $a_{i}=\<f,u_{i}\>$  then 
  \begin{align*}
    I_{1}   &=-\bigg(\sum_{i=1}^{\infty}\lambda_{i}^{3}|a_{i}|^{2}\bigg)\bigg(\sum_{i=1}^{\infty}|a_{i}|^{2}\bigg)+\bigg(\sum\lambda_{i}^{2}|a_{i}|^{2}\bigg)\bigg(\sum \lambda_{i}|a_{i}|^{2}\bigg).
  \end{align*}
  We multiply by
  \[\bigg(\sum \lambda_{i}^{2}|a_{i}|^{2}\bigg)\bigg(\sum \lambda_{i}|a_{i}|^{2}\bigg).\]
  The coefficients $|a_{i}|^{2}$ give us weights for a sequence space $\ell^{2}$, which contains the sequences $\lambda_{i}^{3/2}$, $\lambda_{i}^{1/2}$.  Consequently
  \begin{align*}
  \bigg(\sum \lambda_{i}^{2}|a_{i}|^{2}\bigg)\bigg(\sum \lambda_{i}|a_{i}|^{2}\bigg)  I_{1}
  &=-\|\lambda_{i}^{3/2}\|_{(\ell^{2},|a_{i}|^{2})}^{2}\|\lambda_{i}^{1}\|_{(\ell^{2},|a_{i}|^{2})}^{2}\|\lambda_{i}^{1/2}\|_{(\ell^{2},|a_{i}|^{2})}^{2}\|1\|_{(\ell^{2},|a_{i}|^{2})}^{2}\\
  &\quad+\<\lambda_{i}^{3/2},\lambda_{i}^{1/2}\>_{(\ell^{2},|a_{i}|^{2})}^{2}\<\lambda_{i},1\>_{(\ell^{2},|a_{i}^{2})}^{2}\\
   &\leq 0
  \end{align*}
  by the Cauchy--Schwarz inequality for the $(\ell^{2},|a_{i}|^{2})$ inner product.
  
  We now show that $I_{2}\leq C_{1}(\<B_{t}f_{t},B_{t}f_{t}\>+\<f_{t},f_{t}\>)$. The bound for all of the terms save $\<B_{t}X_{t}f_{t},B_{t}f_{t}\>$ follow directly from the $H^{1}$ boundedness assumptions.  
  \begin{align*}
    \<B_{t}X_{t}f_{t},B_{t}f_{t}\>&=\<[B_{t},X_{t}]f_{t},B_{t}f_{t}\>_{t}+\<X_{t}B_{t}f_{t},B_{t}f_{t}\>_{t}\\
    &=\<[B_{t},X_{t}]f_{t},B_{t}f_{t}\>_{t}-\<\div_{t}X_{t} B_{t}f_{t},B_{t}f_{t}\>.
  \end{align*}
  Consequently
  \begin{align*}
    |\<B_{t}X_{t}f_{t},B_{t}f_{t}\>|&\leq \|[B_{t},X_{t}]\|\|f_{t}\|_{1}\|B_{t}f_{t}\|_{0}+\|\div_{t}X_{t}\|\|f_{t}\|_{1}^{2}\\
    &\leq C(\<B_{t},B_{t}\>+\<f_{t},f_{t}\>).
\end{align*}

Lastly  the term $I_{3}\leq C_{2}\<f_{t},f_{t}\>$, but this is straightforward.

Now we have that
\[\D_{t}\lambda(t,f_{0})\leq C_{1}\lambda(t,f_{0})+ C_{2},\]
so by Gr\"onwall's inequality 
\[\lambda(t,f_{0})\leq e^{C_{1}t}\left[\frac{\<B_{0}f_{0},B_{0}f_{0}\>_{0}}{\<f_{0},f_{0}\>_{0}}+\frac{C_{2}}{C_{1}}\right]-\frac{C_{2}}{C_{1}}.\]
\endproof
\begin{lem}Let $B_{t}$, $X_{t}$ $S_{t}$,$\<\cdot,\cdot\>_{t}$ satisfy the conditions of Lemma \ref{lem:osc_evolution}.  Suppose $\<B_{0}f_{0},B_{0}f_{0}\>\leq \lambda \<f_{0},f_{0}\>$.  Let $f_{t}$ satisfy equation \eqref{eqn:generic_evolution_equation}.  Then
  \[\|f_{t}\|_{1/2}^{2}\geq e^{-(C_{2}\lambda+C_{3})t}\|f_{0}\|_{1}.\]
  \label{lem:evolution_lower_bound}
\end{lem}
\proof
Consider 
\begin{align*}
  \D_{t}[\<f_{t},f_{t}\>+\<B_{t}f_{t},f_{t}\>_{t}]&=\<\dot{B}_{t}f_{t},f_{t}\>+\<B_{t}f_{t},\gamma_{t}f_{t}\>-2\<B_{t}f_{t},B_{t}f_{t}\>-2<B_{t}f_{t},f_{t}\>\\ 
  &\quad+2\<B_{t}X_{t}f_{t},f_{t}\> +2\<B_{t}S_{t}f_{t},f_{t}\>.
  \intertext{Consequently}
  \D_{t}(\<f_{t},f_{t}\>+\<B_{t}f_{t},f_{t}\>_{t})&\geq -C\|f_{t}\|_{1/2}-\frac{\<B_{t}f_{t},B_{t}f_{t}\>}{\<f_{t},f_{t}\>}\<f_{t},f_{t}\>\\
  &\geq -C(\<f_{t},f_{t}\>+\<B_{t}f_{t},f_{t}\>)-C_{2}(\lambda(\<f_{t},f_{t}\>+\<B_{t}f_{t},f_{t}\>),
\end{align*}
by applying Lemma \ref{lem:osc_evolution}.  The result then follows by Gr\"onwall's inequality.
\endproof
\begin{lem}
  \label{lem:ext_adjoin_lb}
  Let $\E|_{\Sigma_{1}}^{*}$ denote the adjoint of the map $H^{1/2}(\D\M)\to H^{1/2}(\D\Sigma_{1})$.  There is a $K\geq 0$  and a $\lambda_{0}$ such that for every  $f\in H^{1/2}(\Sigma_{1})$, $\lambda\geq \lambda_{0}$ if $\|f\|_{0}^{2}\leq \lambda \|f\|_{-1/2}^{2}$ then 
  \[\|\E^{*}(f)\|_{-1/2}^{2}\geq e^{-K\lambda}\|f\|_{-1/2}^{2}.\]
\end{lem}
\proof
The proof of this makes use of Lemmas \ref{lem:extension_restriction_adjoint} and \ref{lem:evolution_lower_bound}.  We append a cylinder $(-\infty,0]\times \D\M$ to $\M$.  We let the metric on $(-\infty,0]\times\D\M$, be equal to a product metric on $\D\M$ after smoothly extending the metric on $\M$.  Denote this new manifold $\tilde{\M}$. The metric tensor on the whole manifold is thus Lipschitz.  Now let $\tilde{\E}$ denote the harmonic extension  from $\D\Sigma_{1}$ to $\tilde{\M}\setminus \Sigma_{1}$.  By similar reasoning to Lemma \ref{lem:extension_comparison}, we have that

\[\D_{\nu\D\M}\E^{0,1}(u)=-\Lambda^{1,0}(\tilde{\E}(u)|_{\D\M})-\Delta^{1/2}_{0}(\tilde{\E}|_{\D\M}),\]
for any $u$.  Setting $u=(\Lambda^{0,1}+\Lambda_{\Sigma_{1}})^{-1}f$, we arrive at
\begin{align*}
  \|\E^{*}(f)\|_{-1/2}&=\|(1+\Delta)^{-1/2}(\Lambda^{0,1}+\Delta^{1/2})(\tilde{\E}|_{\D\M}(\Lambda^{0,1}+\Lambda_{\Sigma_{1}})^{-1}(f))\|_{1/2}^{2}.
\end{align*}
\begin{align*}
  \|(\Lambda^{0,1}+\Lambda_{\Sigma_{1}})^{-1}f\|_{1}^{2}\leq C\lambda\|(\Lambda^{0,1}+\Lambda_{\Sigma_{1}})f\|_{0}^{2}.
\end{align*}
  Hence we can apply Lemma \ref{lem:evolution_lower_bound} to $\tilde{\E}$ to yield that
  \[\|\tilde{\E}((\Lambda^{0,1}+\Lambda_{\Sigma_{1}})^{-1}f)|_{\D\M}\|_{1/2}^{2}\geq e^{-CK\lambda}\|f\|_{-1/2}^{2}.\]
  Lastly we note that $(1+\Delta)^{-1/2}(\Lambda^{1,0}+\Delta^{1/2})$ is a bounded invertible map from $H^{1/2}\to H^{1/2}$.  Giving us the result.
\endproof
\begin{lem}
  \label{lem:iteration_step_existence}
  There are numbers $C(\Sigma_{1},g)$ and $K(\Sigma_{1},g)$ such that  for every $v\in H^{1}(\D\Sigma_{1})$,  $\|v\|_{1/2}=1$, $\|v\|_{1}\leq \lambda$, there is a $w\in H^{1/2}(\D\M)$ such that 
  \begin{equation}
    \|v-\E_{1}(w)|_{\Sigma_{1}}\|^{2}_{1/2}=(1-e^{-K\lambda})
    \label{eqn:approximation_step}
  \end{equation}
  and
  \begin{equation}
    \|v-\E_{1}(w)|_{\Sigma_{1}}\|^{2}_{1}\leq\lambda(1+ Ce^{-K\lambda}).
    \label{eqn:smoothness_step}
  \end{equation}
\end{lem}
\proof Let $A=(I+\Lambda_{1})$
\begin{align*}
  \<A (v-\E(u)|_{\Sigma^{1}}),v-\E(u)\>&=\<A v,A v\>+\<\E(u)|_{\Sigma^{1}},A \E(u)|_{\Sigma_{1}}\>-2\<v,A \E(u)|_{\Sigma_{1}}\>.
\end{align*}
Let $\sigma\in[-1,1]$, then set $u=\sigma  (I+\Lambda_{0})^{-1}(\E|_{\Sigma^{1}})^{*}(A v)$.  Hence
\begin{align*}
  \<\E(u)|_{\Sigma_{1}},Av\>&=\sigma\<(I+\Lambda_{0})^{-1}\E^{*}A(v),\E^{*}A(v)\>\geq \sigma e^{-K\lambda}\|v\|_{1/2}^{2}.
\end{align*}
Whereas $\<\E(u)|_{\Sigma_{1}},A\E(u)|_{\Sigma_{1}}\>_{1}\leq C\sigma^{2}$, independent of the choice of $v$.  We note that $\lambda\geq \lambda_{0}>0$, hence by making $K$ bigger we can make  $C^{-1}\geq e^{-K\lambda}$.  Then setting $\sigma=e^{-K\lambda}$, yields \eqref{eqn:approximation_step}.

For \eqref{eqn:smoothness_step}, let $w=\E(u)|_{\D\Sigma}$. 
\begin{align*}
  \frac{\<A (v-w),A (v-w)\>}{\<A (v-w),v-w\>} \hspace{-8em}&\\
 &\leq \frac{\<A v,v\>}{\<A (v-w),v-w\>}\Bigg[\frac{\<A v,A v\>}{\<A v,v\>}+2\frac{\<A v,A v\>^{1/2}\<A w,A w\>^{1/2}}{\<A v,v\>}\\
 &\quad+\frac{\<A w,A w\>}{\<A v,v\>}\Bigg]\\
  &\leq (1-\sigma)\Bigg[\lambda+2\lambda^{1/2}\Bigg(\frac{\<A w,A w\>}{\<A w,w\>}\frac{\<A w,w\>}{\<A v,v\>}\Bigg)^{1/2}+\frac{\<A w,A w\>}{\<A w,w\>}\frac{\<A w,w\>}{\<A v,v\>}\Bigg].
\end{align*}
By Lemma \ref{lem:osc_evolution} we know that 
\[\frac{\<A w,A w\>}{\<A w,w\>}\leq C \lambda,\]
And we know that $\<A w,w\>\leq (Ce^{-K\lambda})^{2}$, because of our choice of $\sigma$. Hence 
\[\frac{\<A (v-w),A (v-w)\>}{\<A (v-w),v-w)\>}\leq \lambda(1+Ce^{-K\lambda}).\]
\endproof

\section{The proofs of the main theorems}
\proof[Proof of Theorem 2]
The idea is a simple iteration scheme applying Lemma \ref{lem:iteration_step_existence}.  Let $\|v_{i}\|_{1/2}=1$.  We construct $u_{i}$ such that $\|v_{i}-\E^{1}(u_{i})|_{\D\Sigma_{1}}\|_{1/2}=(1-\mu_{i})^{1/2}$  Then we set $v_{i+1}=(v_{i}-\E^{1}(u_{i}))/(1-\mu_{i})^{1/2}$. Then 
\[\frac{\|v_{i+1}\|_{1}^{2}}{\|v_{i+1}\|_{1/2}^{2}}\leq \lambda_{i+1},\]
Where
\[\lambda_{i+1}=\lambda_{i}(1+Ce^{-K\lambda_{i}}),\]
and $\mu_{i}=e^{-K\lambda_{i}}$.  Then at step $i$ we consider
\[\|v-\E(\tilde{u}_{i})\|_{1/2}^{2}=\prod_{k=1}^{i}(1-\mu_{k}),\]
where 
\[\tilde{u}_{i}=\sum_{k=1}^{i}\prod_{j=1}^{k-1}(1-\mu_{j})^{1/2}u_{k}.\]
Recall that $\|u_{i}\|\leq C\|v_{i}\|\leq C$, so 
\[\|\tilde{u}_{i}\|\leq C\sum_{k=1}^{i}\prod_{j=1}^{k-1}(1-\mu_{j})^{1/2}.\]
First we show that 
\[\prod_{k=1}^{i}(1-\mu_{k})\] 
tends to zero, as $i$ tends to infinity.  To do this, consider the logarithm
\[\log\left(\prod_{k=1}^{i}(1-\mu_{k})\right)\leq \sum_{k=1}^{i}\log(1-\mu_{k})\leq -\sum_{k=1}^{i}\mu_{k}\]

Now we need to determine lower bounds for $\mu_{k}$.  We examine
\[K\lambda_{k+1}=K\lambda_{k}(1+Ce^{-K\lambda_{k}}).\]
Setting $\sigma_{k}=K\lambda_{k}$, we can observe that this is the Euler scheme for approximating the differential equation
\[f'(t)=Cf(t)e^{-f}.\]
This has as its solution
\[f(t)=\log(\li^{-1}(Ct+D)),\]
where $\li$ is the logarithmic integral
\[\li(t)=\int_{2}^{t}\frac{1}{\log(\tau)}\:d\tau.\]
\[(\li^{-1})'(t)=\frac{1}{\li'(\li^{-1}(t))}=\log(\li^{-1}(t)).\]

So we claim 
\[\log(\li^{-1}(Ck+\li(\exp(\sigma_{0}))))\leq\sigma_{k}\leq \log(\li^{-1}(Ck+\li(\exp(\sigma_{0}))))+1\]
And the proof is by induction.  Let $\theta=\log(\li^{-1}(w+D))$. Assume $\theta e^{-\theta}\leq (12C)^{-1}$, let $h=y-\log(\li^{-1}(w+D))\geq $, let $y'=y+Cye^{-y}$, and let $h'=y'-\log(\li^{-1}(w+C+D))$. Because $\log(\li^{-1}(\cdot))$ is a concave function
\[\log(\li^{-1}(w+C+D))\leq \log(\li^{-1}(w+D))+\frac{\log(\li^{-1}(w+D))}{\li^{-1}(w+D)}C\leq \theta+C\theta e^{-\theta}.\]
Then 
\begin{align*}
  h'&\geq h+Ce^{-y}y-C\theta e^{-\theta}\\
  &=h+Ce^{-y}(y-\theta)+C(e^{-(y-\theta)}-1)e^{-\theta}\theta\\
  &\geq h(1+Ce^{-y}-Ce^{-\theta}\theta)\\
  &\geq h/2\geq 0.
\end{align*}
We let $w=CK$ and $y=\sigma_{k}$, it follows that 
  \[\sigma_{k+1}\geq\log(\li^{-1}(C(k+1)+D))\]
For the upper bound, we once again use the fact that $\log(\li^{-1}(t))$ is concave, and that
\[(\log\circ\li^{-1})''(t)=\frac{1}{\li^{-1}(t)^{2}}\left[\log(\li^{-1}(t))-\log(\li^{-1}(t))^{2}\right]\]
Consequently, 
\begin{align*}
  \log(\li^{-1}(w+C+D)&\geq\theta+C\theta e^{-\theta}-\theta(\theta-1)e^{-2\theta}C^{2}.
\end{align*}
So
\begin{align*}
  h'\leq h+Che^{-y}-C(1-e^{-h})\theta e^{-\theta}+\theta(\theta-1)e^{-2\theta}C^{2}.
\end{align*}
If $h\geq 3\theta e^{-\theta} C$, then 
\[h'\leq (h+C e^{-\theta} h-C h \theta e^{-\theta}/2+Ch\theta e^{-\theta}/3\leq h(1+C(e^{-\theta}(1-\theta/6))<h\leq 1. \]
If $h\leq 3\theta e^{-\theta} C$, then
\[h'\leq 1/4+1/4+1/8+1/4<1.\]
So once again if $\sigma_{k}=y$ then $\sigma_{k+1}=y'$ and if
\[\sigma_{k}\leq \log(\li^{-1}(Ck+D))+1\],
then 
\[\sigma_{k+1}\leq \log(\li^{-1}(C(k+1)+D))+1.\]
Consequently, to get towards the initial induction step, set 
\[\sigma_{0}=\max\{K\|v\|_{1}^{2}/\|v\|_{1/2},(e^{x}/x)^{-1}(12C))\}.\]
Because $D=\li(\exp(\sigma_{0}))$, we have that $y_{0}=\sigma_{0}$, and hence by induction the upper and lower bounds follow.

Now $\prod_{k=1}^{i}(1-\mu_{k})\leq \exp(-\sum_{k=1}^{i}(\li^{-1}(Ck+D))^{-1})$, so
\begin{align*}
  \prod_{k=1}^{i}(1-\mu_{k})&\leq \exp(-\int_{0}^{i}\frac{1}{\li^{-1}(Ct+D)}\:dt\\
  &\leq \exp\left( -C \int_{\li^{-1}(D)}^{\li^{-1}(Ci+D)}\frac{ds}{s\log(s)}\right)\\
  &\leq \exp(-[\log\log(\li^{-1}(Ci+D)-\log(\sigma_{0})]\\
  &\leq \sigma_{0}\frac{1}{\log(\li^{-1}(Ci+D))},
\end{align*}
by setting $s=\li^{-1}(Ct+D)$, the $ds=C\log(s)dt$.  Now for $\varepsilon>0$  let
\[i=\li(\exp(( (\sigma_{0}/\varepsilon))-\li(\exp(\sigma_{0})))/C\]
And consequently
\[\|\tilde{u}_{i}\|^{2}\leq C\sum_{k=1}\frac{\sigma_{0}}{\log(\li^{-1}(Ck+D))}\leq\int_{\sigma_{0}}^{(\sigma_{0}/\varepsilon)}(e^{t}/t)\leq e^{(\sigma_{0}/\varepsilon)^{\alpha}-2\log(\sigma_{0})}.\]

\endproof
\begin{lem}
  \label{lem:boundary_inversion}
 t Let $\M$ and $\Sigma\subset \M$  be compact $C^{2}$-smooth manifolds with boundary. There exist numbers and $r_{0},c>0$ such that for every $x\in \D\Sigma$ there is a map $\Phi:\D\M\times[0,1]\to \M$ with $g(\D_{t}\varphi,\nu)\geq c$, and $\D\Sigma_{1}\supset \supset B(x,r_{0})\cap \D\Sigma$.
\end{lem}
\proof
By the collar neighbourhood theorem we have a map $\varphi:\D\M\times[0,1]\to \M\setminus \Sigma$ which is diffeomorphic onto its image.  Now fix a point $x_{0}$ in $\D\M$ and extend the path $t\mapsto \varphi(x_{0},t)$ to a path $\gamma:[0,2]\to \M$ from $x_{0}$ to $x\in \D\Sigma$.  By the tubular neighbourhood theorem their is a tubular neighbourhood of $\gamma$ in $\Sigma_{1}\setminus \Sigma$ diffeomorphic to $[0,1]\times B$ where $B\subset \D\Sigma_{1}$ is a neighbourhood of $\varphi(x_{0},1)$ via a map $\psi$.  The diffeomorphism can be choosen so that for $x$ in $B'=(x'\mapsto (\varphi(x',1))^{-1}(B)$ the map 
\[\Psi:(x,t)\mapsto \begin{cases}
    \varphi(x,t) & t\in [0,1]\\
    \psi(\varphi(x,1),t-1) & t\in[0,2]
\end{cases},\]
is a $C^{2}$ diffeomorphism of $B'\times [0,2]$ onto its image.  Let $\zeta:B_{n-1}(0,\varepsilon)\to B'$ be a $C_{2}$ diffeomorphism.  Then we consider a monotone decreasing map $\rho:[0,\varepsilon]\to \RR+$ equal to $1$ in a neighbourhood of $0$ and $0$ in a neighbourhood of $\varepsilon$.  With this we define
\[\Phi(x,t)=\begin{cases}
  \varphi(x,t) & x\in \D\M\setminus B'\\
  \psi(x,2t\rho(|\zeta^{-1}(x)| +(1-\rho(|\zeta^{-1}(x)|)) t), & x \in B'
\end{cases}.\]
\endproof

\endproof

\begin{lem}
  \label{lem:inner_boundary_functions}. There is a number $C=C(\M,g,\Sigma)\geq 0$ such that for every piecewise continuous metric tensor $h$ on $\Sigma$ and every $x\in B(x',r_{0}/2)$ there are Lipschitz functions $u_{1}$ and $v_{1}$ supported in $B(x',r_{0})$ satisfying $\|u_{1}\|_{1/2}=1$ and $\|v_{1}\|_{1/2}=1$, such that 
  \[(\Lambda_{\Sigma_{1},g}-\Lambda_{\Sigma_{1},h})(u_{1})(v_{1})\geq C\frac{|h(x)-g(x)|}{|g(x)|}.\]
\end{lem}
\proof
%\begin{align*}
%  \int_{\M}\xi^{2} r_{0}^{(n-1)/2}
%\end{align*}
The idea for this proof was adapted from the $p$-harmonic case \cite{p-harmonic_Calderon}. Consider $d(\psi u)\wedge \star_{g} dv$ for $v\in H^{1}_{0}$, where $\psi$ is a smooth cutoff, and $u\circ \varphi^{-1}=|\xi|^{-1/2}e^{\xi\cdot x}e^{-|\xi|\tau}$ where $\varphi$ is some coordinate map.   Suppose $\psi$ is supported in $B(0,r_{0})$ in this coordinate map.  Let $\star$ be the pullback of the Euclidean Hodge star under $\varphi$, and assume $\star_{g}=\star$ at $\varphi^{-1}(0)$, then 
\begin{align*}
  \bigg| \int_{\Sigma^{1}}d(\psi u)\wedge\star_{g} dv\bigg| &= \bigg|\int_{\Sigma^{1}}du\wedge\star d(\psi v)+du\wedge(\star-\star_{g}) \psi dv+u d\psi\wedge\star_{g}dv- u\wedge d\star vd\psi\bigg |\\
&\leq C(\|u\|_{H^{1}(B(0,r_{0}))}r_{0}+1/(|\xi|r_{0})+1/(|\xi|^{1/2}r_{0}^{n/2-1}))\|v\|_{1}.
\end{align*}
while $\|u\psi\|_{1}\sim r_{0}^{(n-1)/2}$.
Consider the orthogonal decomposition $H^{1}_{0}(\M)$ to $H^{1}_{0}(\varphi^{-1}(B(0,r_{0})\times[0,r_{0}]))\oplus \A$ where $\A$ is the space of functions which are zero on $\D\M$ and harmonic in $\varphi^{-1}(B_{r_{0}}\times[0,r_{0})$.
  So by making $|\xi|^{1/2}r_{0}^{n-3/2}$ big  and $r_{0}$ small, we can make $\Delta_{g} \psi u$ arbitrarily small in $H^{-1}$ relative to the $H^{1}$ norm of $\psi u$.  So the $H^{1}$ difference between $\psi u$ and $\E(\psi u)$ is also arbtitrarily small.  Hence if we divide $\psi u$ by $r_{0}^{(n-1)/2}$ we get
  \[\int_{\M}d\psi u\wedge\star_{g}d(\psi u)\to 1 \]
as $|\xi|^{1/2}r_{0}^{n-3/2}\to \infty$ and $r_{0}\to 0$.
and hence 
\[\int_{\D\M}\Lambda_{g}(\psi u)(\psi u)\:d\vol_{g}\to 1.\]
By homogeneity we get that  
\[\int_{\D\M}\Lambda_{h}(\psi u)(\psi u)\:d\vol_{h}\to |\xi|_{h}/|\xi|_{g},\]
where the rate of convergence is dependent only on the Lipschitz constants of the metric tensors in a neighbourhood of $x$.

%The $H^{1}$ norm of $u\psi$ is asymptotically given by $|\xi|^{1/2}$,  $r_{0}=\varepsilon/3C$,  then $\xi^{1/2}\geq 3C\varepsilon^{-2}$.  Then $\epsilon<|1-|\xi|_{2}/|\xi|_{1}|/3$.
\endproof

\proof[Proof of Theorem 1]
We cover $\D\Sigma$ in balls of sufficiently small radius $B(x_{i},r_{i})$ such that we can apply Lemma \ref{lem:boundary_inversion} on the set $\Sigma\cap B(x_{i},2r_{i})$. By compactness there is a finite such collection.  Suppose $|g-h|_{g}=|1-h(X,X)/g(X,X)|$ for $X\in T_{x}\M$ achieves its maximum  in $\D\Sigma$ at $x\in B(x_{i},r_{i})$. Take $u_{1}$ and $u_{2}$ given by Lemma \ref{lem:inner_boundary_functions}.  As per Lemma \ref{lem:extension_comparison}, take $v=(\Lambda_{\Sigma_{1},g}+\Lambda^{0,1}_{g})^{-1}(\Lambda_{\Sigma_{1},h}+\Lambda^{0,1}_{\M\setminus \Sigma_{1}})u_{2}$. Let $\varepsilon=\|h(x)-g(x)\|/C''$, where $C''=C''(K,R,g,\Sigma_{1})$ will be determined later.  Then let $u^{\varepsilon}$ $v^{\varepsilon}$ be  Lipschitz approximations of $u$ and $v$ in  $H^{1/2}$ satisfying
\[\|v-v^{\varepsilon}\|_{1/2}\leq \varepsilon,\]
and
\[\|v^{\varepsilon}\|_{1,\infty}\leq C'/\varepsilon.\]

Then there are a $u_{0}$ $v_{0}$ such that $\|\E^{1}(u_{0})-u_{\varepsilon}\|_{1/2}\leq \varepsilon $ given by Theorem \ref{thm:harmonic_density_estimate}.  Consequently 
\begin{align*}
  (\Lambda^{1}&-\Lambda^{2})(u_{0})(v_{0})\\
  &=(\Lambda^{1}_{\Sigma_{1}}-\Lambda^{2}_{\Sigma_{1}})(\E^{1}(u_{0})|_{\Sigma_{1}})(\Lambda_{\Sigma_{1},g}+\Lambda^{0,1}_{g})^{-1}(\Lambda_{\Sigma_{1},h}+\Lambda^{0,1}_{\M\setminus\Sigma_{1}})\E^{1}(v_{0})|_{\Sigma_{1}}\\
  &\geq (\Lambda^{1}_{\Sigma_{1}}-\Lambda^{2}_{\Sigma_{1}})(u_{1})(u_{2})- (\|\Lambda^{1}_{\Sigma_{1}}\|+\|\Lambda^{2}_{\Sigma_{1}}\|)\|\E_{1}u_{0}|_{\Sigma_{1}}-u_{1}\|\E_{2}v_{0}|_{\Sigma_{1}}\|\\
  &\quad-\|u_{1}\|\|A\|\|\E_{1}(v_{0})-v\|\\
  &\geq C(g,\Sigma_{1},K)|h-g|_{g}-C'(R,g,\Sigma_{1},K)\varepsilon\\
  &\geq \frac{\|h-g|}{C(g,K,\Sigma_{1})},
\end{align*}
where we have chosen $C''(g,K,R,\Sigma_{1})$, such that $C'/C''<C/2$ 
Dividing through by the upper bounds on $u_{0}$ and $v_{0}$ yields the result, 
\[\|\Lambda_{1}-\Lambda_{2}\|\geq C_{1}\|h-g\|_{g}\exp(2C_{2}\|h-g\|_{g}^{-2C^{3}}).\]
 We choose the worst constants over all $i$ assuming  $x\in B(x_{i},r_{i})$.    Then we take the logarithm of both sides yielding
 \[|\log(\|\Lambda_{1}-\Lambda_{2}\|)|\leq |\log C_{1}\|h-g\|_{g}|+2C_{2}\|h-g\|_{g}^{-2C_{3}}.\]
 Lastly $x^{-\beta}$ dominates $|\log(x)|$ for  every $\beta$ and every $x\leq R$. We subsume the $|\log C_{1}\|h-g\|_{g}|$ into the $\|h-g\|_{g}^{-2C^{3}}$ with an appropriate constant (dependent on $R$), to yield the main theorem.
\endproof

\bibliographystyle{plain}
%\bibliography{../Jan}

\begin{thebibliography}{10}

\bibitem{Calderon_in_plane}
Kari Astala and Lassi P{\"a}iv{\"a}rinta.
\newblock Calder\'on's inverse conductivity problem in the plane.
\newblock {\em Ann. of Math. (2)}, 163(1):265--299, 2006.

\bibitem{anisotropic_Calderon_plane}
Kari Astala, Lassi P{\"a}iv{\"a}rinta, and Matti Lassas.
\newblock Calder\'on's inverse problem for anisotropic conductivity in the
  plane.
\newblock {\em Comm. Partial Differential Equations}, 30(1-3):207--224, 2005.

\bibitem{Collar_neighbourhood}
Morton Brown.
\newblock Locally flat imbeddings of topological manifolds.
\newblock {\em Ann. of Math. (2)}, 75:331--341, 1962.

\bibitem{Calderon_Problema_classic}
Alberto-P. Calder{\'o}n.
\newblock On an inverse boundary value problem.
\newblock In {\em Seminar on {N}umerical {A}nalysis and its {A}pplications to
  {C}ontinuum {P}hysics ({R}io de {J}aneiro, 1980)}, pages 65--73. Soc. Brasil.
  Mat., Rio de Janeiro, 1980.

\bibitem{Calderon_Problem_new}
Alberto~P. Calder{\'o}n.
\newblock On an inverse boundary value problem.
\newblock {\em Comput. Appl. Math.}, 25(2-3):133--138, 2006.

\bibitem{CaroRogers_Lipschitz_conductivities}
Pedro Caro and Keith Rogers.
\newblock Global uniqueness for the {C}alder\'on problem with lipschitz
  conductivities.
\newblock {\em arXiv}.

\bibitem{Limiting_Carleman_weights}
David Dos Santos~Ferreira, Carlos~E. Kenig, Mikko Salo, and Gunther Uhlmann.
\newblock Limiting {C}arleman weights and anisotropic inverse problems.
\newblock {\em Invent. Math.}, 178(1):119--171, 2009.

\bibitem{skeletal_muscle_anisotropy}
B.R. Epstein and K.R. Foster.
\newblock Anisotropy in the dielectric properties of skeletal muscle.
\newblock {\em Medical and Biological Engineering and Computing}, 21(1):51--55,
  1983.

\bibitem{GilbargTrudinger}
David Gilbarg and Neil~S. Trudinger.
\newblock {\em Elliptic partial differential equations of second order}.
\newblock Classics in Mathematics. Springer-Verlag, Berlin, 2001.
\newblock Reprint of the 1998 edition.

\bibitem{Lipschitz_conductivities}
Boaz Haberman and Daniel Tataru.
\newblock Uniqueness in {C}alder\'on's problem with {L}ipschitz conductivities.
\newblock {\em Duke Math. J.}, 162(3):496--516, 2013.

\bibitem{cat_dorsal_column_anisotropy}
James B.~Ranck Jr. and Spencer~L. BeMent.
\newblock The specific impedance of the dorsal columns of cat: An anisotropic
  medium.
\newblock {\em Experimental Neurology}, 11(4):451 -- 463, 1965.

\bibitem{Lassas_Uhlmann_Taylor}
Matti Lassas, Michael Taylor, and Gunther Uhlmann.
\newblock The {D}irichlet-to-{N}eumann map for complete {R}iemannian manifolds
  with boundary.
\newblock {\em Comm. Anal. Geom.}, 11(2):207--221, 2003.

\bibitem{Lee_smooth_manifolds}
John~M. Lee.
\newblock {\em Introduction to smooth manifolds}, volume 218 of {\em Graduate
  Texts in Mathematics}.
\newblock Springer, New York, second edition, 2013.

\bibitem{p-harmonic_Calderon}
Mikko Salo and Xiao Zhong.
\newblock An inverse problem for the {$p$}-{L}aplacian: boundary determination.
\newblock {\em SIAM J. Math. Anal.}, 44(4):2474--2495, 2012.

\bibitem{Sylvester_layer_stripping}
John Sylvester.
\newblock A convergent layer stripping algorithm for the radially symmetric
  impedance tomography problem.
\newblock {\em Comm. Partial Differential Equations}, 17(11-12):1955--1994,
  1992.

\bibitem{Sylvester_Uhlmann}
John Sylvester and Gunther Uhlmann.
\newblock A global uniqueness theorem for an inverse boundary value problem.
\newblock {\em Ann. of Math. (2)}, 125(1):153--169, 1987.

\bibitem{psidos_NLPDE}
Michael~E. Taylor.
\newblock {\em Pseudodifferential operators and nonlinear {PDE}}, volume 100 of
  {\em Progress in Mathematics}.
\newblock Birkh\"auser Boston, Inc., Boston, MA, 1991.

\bibitem{nonlinear_fnl_analysis_and_apps_IIa}
Eberhard Zeidler.
\newblock {\em Nonlinear functional analysis and its applications. {II}/{A}}.
\newblock Springer-Verlag, New York, 1990.
\newblock Linear monotone operators, Translated from the German by the author
  and Leo F. Boron.

\end{thebibliography}
\def\cprime{$'$}

\end{document}